\documentclass{amsart}

\usepackage{txfonts}
\newtheorem{theorem}{Theorem}[section]
\newtheorem{proposition}[theorem]{Proposition}
\newtheorem{lemma}[theorem]{Lemma}

\newtheorem{corollary}[theorem]{Corollary}

\theoremstyle{definition}
\newtheorem{definition}[theorem]{Definition}

\theoremstyle{remark}
\newtheorem{remark}[theorem]{Remark}

\numberwithin{equation}{section}

\catcode`\ç=13
\defç{\c{c}}
\catcode`\é=13
\defé{\'e}
\catcode`\à=13
\defà{\`a}
\catcode`\è=13
\defè{\`e}
\catcode`\â=13
\defâ{\^a}
\catcode`\ù=13
\defù{\`u}
\catcode`\ê=13
\defê{\^e}
\catcode`\î=13
\defî{\^\i}
\catcode`\ô=13
\defô{\^o}
\newcommand{\field}[1]{\mathbb{#1}}

\newcommand{\Q }{\field{Q}}
\newcommand{\Z }{\field{Z}}

\def\wdim{{\rm wgldim}}
\def\lgldim{{l\rm.gldim}}
\def\rgldim{{r\rm.gldim}}

\def\lGwdim{{l\rm.wGgldim}}
\def\lGgldim{{l\rm.Ggldim}}

\def\rGgldim{{r\rm.Ggldim}}

\def\lFPD{{ l\rm.FPD}}

\def\pd{{\rm pd}}
\def\fd{{\rm fd}}
\def\id{{\rm id}}

\def\Gpd{{\rm Gpd}}
\def\Gfd{{\rm Gfd}}
\def\Gid{{\rm Gid}}

\def\Ker{{\rm Ker}}

\def\Ext{{\rm Ext}}

\def\Hom{{\rm Hom}}

\def\sup{{\rm sup}}
\begin{document}

\title{Global Gorenstein Dimensions}

\author[D. Bennis]{Driss Bennis}
\address{Department of Mathematics, Faculty of Science and Technology of Fez, Box 2202, University S. M.
Ben Abdellah Fez, Morocco}
 \email{driss\_bennis@hotmail.com}

\author[N. Mahdou]{Najib Mahdou}
\address{Department of Mathematics, Faculty of Science and Technology of Fez, Box 2202, University S. M.
Ben Abdellah Fez, Morocco}
 \email{mahdou@hotmail.com}



\dedicatory{Dedicated to our Advisor Salah-Eddine Kabbaj.}


\begin{abstract} In this paper, we prove that the global Gorenstein projective dimension of a ring
$R$ is equal to the global Gorenstein injective dimension of $R$,
and that the  global Gorenstein flat dimension  of $R$ is smaller
than the common value of the terms of this equality.
\end{abstract}

\footnotetext[0]{ *Proceedings of the American Mathematical
Society, To appear }\value{counter}\setcounter{footnote}{0}

\maketitle

\section{Introduction} Throughout this paper, $R$ denotes
a non-trivial associative ring with identity, and all modules are,
if not specified otherwise, left $R$-modules. All the results,
except Proposition \ref{pro-QFrings}, are formulated for left
modules and the corresponding results for right modules hold as
well. For an $R$-module $M$, we use $\pd_R(M)$, $\id_R(M)$, and
$\fd_R(M)$ to denote, respectively, the classical projective,
injective, and flat dimension  of $M$.  We use $\lgldim(R)$ and
$\rgldim(R)$ to denote, respectively, the classical left and right
global dimension of $R$, and   $\wdim(R)$ to denote the weak
global dimension  of $R$. Recall that the left finitistic
projective dimension of R is the quantity
$\lFPD(R)=\sup\{\pd_R(M)\;|\;M\;is\;an \;
 R\!-\!module\;with$ $\;\pd_R(M)<\infty\} $.\\ Furthermore, we use $\Gpd_R(M)$,
$\Gid_R(M)$, and $\Gfd_R(M)$ to denote, respectively, the
Gorenstein projective, injective, and flat dimension  of $M$ (see
 \cite{LW, Rel-hom, HPhD}).\bigskip

The main result of this paper is an analog of a classical equality
that is used to define the global dimension of $R$, see
\cite[Theorems  9.10]{Rot}. For Noetherian rings the following
theorem is proved in \cite[Theorem 12.3.1]{Rel-hom}.

\begin{theorem}\label{thm-main} The following equality holds:
  $$\sup\{\Gpd_R(M)\,|\,M\;is\;an\;R\!-\!module\}  =
\sup\{\Gid_R(M)\,|\,M\;is\;an\;R\!-\!module\}.$$
\end{theorem}
We call the common value of the quantities in the theorem the
\textit{left Gorenstein global dimension} of $R$ and denote it by
$\lGgldim(R)$. Similarly, we set
$$\lGwdim(R)=\sup\{\Gfd_R(M)\,|\,M\;is\;an\;R\!-\!module\}$$ and call this quantity the \textit{left weak Gorenstein global dimension} of
$R$.

\begin{corollary}\label{cor-main} The following inequalities hold:
\begin{enumerate}
        \item $\lGwdim(R)\leq \sup\{\lGgldim(R),\rGgldim(R)\}. $
    \item  $\lFPD(R)\leq \lGgldim(R)\leq \lgldim(R)$.
      \item   $ \lGwdim(R)\leq  \wdim(R).$
\end{enumerate}
Equalities hold in (2) and (3) if $\wdim(R)<\infty$.
\end{corollary}

The theorem and its corollary are proved in Section 2.

\section{Proofs of the main results}

The proof   use the following results:

\begin{lemma}\label{GPDlem} If $\sup\{\Gpd_R(M)\,|\,M\;is\;an\;R\!-\!module\}<\infty$, then, for  a
positive integer $n$,  the following are equivalent:
\begin{enumerate}
    \item $\sup\{\Gpd_R(M)\,|\,M\;is\;an\;R\!-\!module\}\leq n$;
    \item   $\id_{R}(P )\leq n$ for every
$R$-module  $P $ with finite projective dimension.
\end{enumerate}
\end{lemma}
\begin{proof}  Use \cite[Theorem 2.20]{HH} and \cite[Theorem 9.8]{Rot}.\end{proof}

The proof of the main theorem depends   on the notions of strong
Gorenstein projectivity and injectivity, which were introduced in
\cite{BM} as follows:

\begin{definition}\textnormal{(\cite[Definition
2.1]{BM}).\label{def-SG-pro-inj}
An $R$-module $M$ is called strongly Gorenstein projective, if
there exists an exact sequence of  projective $R$-modules
$$ \mathbf{P}=\
\cdots\stackrel{f}{\longrightarrow}P\stackrel{f}{\longrightarrow}P\stackrel{f}{\longrightarrow}P
\stackrel{f}{\longrightarrow}\cdots
    $$  such that  $M\cong \Ker\, f$ and such that $\Hom_R ( -, Q) $ leaves the
sequence $\mathbf{P}$ exact whenever $Q$ is a projective $R$-module.\\
\indent Strongly Gorenstein injective modules are defined dually.}
\end{definition}

\begin{remark}\label{Cara2-Gfor}
It is easy to see that  an $R$-module $M$ is  strongly Gorenstein
projective if and only if there exists a short exact sequence of
$R$-modules $0\rightarrow M\rightarrow P\rightarrow M\rightarrow
0,$ where $P$ is  projective, and $\Ext^i_R(M,Q)=0$ for some
integer  $i>0$ and for every $R$-module $Q$ with finite
projective dimension (or for every projective $R$-module $Q$).\\
\indent Strongly Gorenstein injective modules are characterized in
similar terms.
\end{remark}

The principal role of these modules is to characterize the
Gorenstein projective and injective modules, as
follows\,\footnote{\; In  \cite{BM} the base ring is assumed to be
commutative. However, for the result needed here, one can show
easily that this assumption is not necessary.}:

\begin{lemma}\textnormal{(\cite[Theorems 2.7]{BM}).}   \label{tcara-G-Gfor}
An $R$-module is Gorenstein projective (resp., injective) if and
only if it is a direct summand of a strongly Gorenstein projective
(resp., injective) $R$-module.
\end{lemma}

\begin{proof}[Proof of Theorem \ref{thm-main}]  For every integer $n$ we need to
show: $$ \Gpd_R(M)\leq n\ for\ every\ R-module\ M\ \quad
\Longleftrightarrow\quad \Gid_R(M)\leq n \ for \ every\ R-module \
M.$$  We only prove the direct implication; the converse one has a
dual proof.\\
Assume first that $M$ is strongly Gorenstein projective. By Remark
\ref{Cara2-Gfor} there is a short exact sequence  $0\rightarrow M
\rightarrow P \rightarrow M \rightarrow 0$ with $P$ is projective.
The Horseshoe Lemma, see [10, Remark page 187], gives a
commutative diagram

$$\begin{array}{ccccccccc}
    &   & 0 &   & 0 &   & 0 &  &   \\
     &   & \downarrow &   & \downarrow  &   & \downarrow &   &  \\
   0 & \rightarrow & M & \rightarrow & P& \rightarrow&M &\rightarrow & 0 \\
     &   & \downarrow &   & \downarrow  &   & \downarrow &   &  \\
   0 & \rightarrow & I_0 & \rightarrow & I_0\oplus I_0 & \rightarrow& I_0 & \rightarrow &0 \\
     &   & \downarrow &   & \downarrow  &   & \downarrow &   &  \\
     &   & \vdots &   & \vdots  &   & \vdots &   &  \\
     &   & \downarrow &   & \downarrow  &   & \downarrow &   &  \\
  0 & \rightarrow & I_n & \rightarrow & E_n& \rightarrow&I_n &\rightarrow & 0 \\
     &   & \downarrow &   & \downarrow  &   & \downarrow &   &  \\
     &   & 0 &   & 0 &   & 0 &  &   \\
  \end{array}
$$
where  $I_i$ is  injective  for $i=0,...,n-1$. Since $P$ is
projective, $\id_R(P)\leq n$ (by Lemma  \ref{GPDlem}), hence $E_n$
is injective. On the other hand, from \cite[Theorem 2.2]{HH1},
$\pd_R(E)\leq n$ for every injective $R$-module $E$. Then,
$\Ext^i_R(E, I_n)=0$ for all $i\geq n+1$. Then, from Remark
\ref{Cara2-Gfor}, $I_n$ is strongly Gorenstein injective, and so
$\Gid_R(M)\leq n$. This implies, from \cite[Proposition 2.19]{HH},
that $\Gid_R(G)\leq n$ for any Gorenstein projective $R$-module
$G$, since every Gorenstein projective $R$-module is a direct
summand of a strongly  Gorenstein projective $R$-module (Lemma \ref{tcara-G-Gfor}).\\
Finally, consider an $R$-module $M$ with $\Gpd_R(M)\leq m\leq n$.
We can assume that $\Gpd_R(M)\not =0$. Then, there exists a short
short exact sequence $0\rightarrow K \rightarrow N \rightarrow M
\rightarrow 0$ such that  $N$ is   Gorenstein projective and
$\Gpd_R(K)\leq m-1$ \cite[Proposition 2.18]{HH}. By induction,
$\Gid_R(K)\leq n$ and $\Gid_R(N)\leq n$. Therefore, using
\cite[Theorems 2.22 and 2.25]{HH} and the long exact sequence of
Ext, we get that $\Gid_R(M)\leq n$.
\end{proof}

\begin{proof}[Proof of Corollary \ref{cor-main}] (1). We may assume that $ \sup\{\lGgldim(R)
, \rGgldim(R)\} < \infty$. Then,  the character module,
$I^*=\Hom_{\Z}(I, \Q/\Z)$, of every injective right $R$-module $I$
has finite projective dimension (by \cite[Theorem 2.2]{HH1} and
\cite[Theorem 3.52]{Rot}). Then, similarly to the proof of
\cite[Proposition 3.4]{HH}, we get that every Gorenstein
projective $R$-module is Gorenstein flat. Therefore,
$\lGwdim(R)\leq
\sup\{\lGgldim(R),\rGgldim(R)\}$.\\[0.2cm]
$(2)$ and $(3)$.  The inequality $\lFPD(R)\leq \lGgldim(R)$
follows
from \cite[Theorem 2.28]{HH}.\\
The inequalities $\lGgldim(R)\leq \lgldim(R)$  and $ \lGwdim(R)\leq
\wdim(R) $ hold true since every projective  (resp.,
flat) module is Gorenstein projective (resp., Gorenstein flat).\\
If $\wdim(R)<\infty$, then, from \cite[Corollary 3]{Jensen},
$\lFPD(R)= \lGgldim(R)= \lgldim(R)$  and, from \cite[Corollary
3.8]{BM}, $ \lGwdim(R)= \wdim(R)$.
\end{proof}

\begin{remark}\label{rem-dif-l.r.Ggldim}
It is well-known that there are examples of rings for which the
left and right global dimensions differ (see  \cite[pages
74-75]{Triv} and \cite{Jat}). Then, by Corollary \ref{cor-main},
the same examples show that there are also examples of rings for
which the left and right Gorenstein global dimensions differ.
However, as the classical case  \cite[Corollary 9.23]{Rot}, we
have $\lGgldim(R)=\rGgldim(R)$ if  $R$ is Noetherian \cite[Theorem
12.3.1]{Rel-hom}.\end{remark}

For the case where $\lGgldim(R)=0$ or $\rGgldim(R)=0$, we have the
following result which is \cite[Theorem 2.2]{BMO} in non-commutative
setting. Recall that a ring is called quasi-Frobenius, if it is
Noetherian and both left and right self-injective (see \cite{NY}).

\begin{proposition}\label{pro-QFrings}
The following are equivalent:
\begin{enumerate}
    \item  $R$ is quasi-Frobenius;
    \item  $\lGgldim(R)=0$;
    \item $\rGgldim(R)=0$.
\end{enumerate}
\end{proposition}
\begin{proof}  The implications  $1\Rightarrow 2$ and $1\Rightarrow3$ are
well-known (see, for example, \cite[Exercise 5, page
257]{Rel-hom}).\\ The implication $2\Rightarrow 1$   follows from
Lemma \ref{GPDlem} and Faith-Walker Theorem \cite[Theorem
7.56]{NY}. The implication  $3\Rightarrow 1$ is proved
similarly.\end{proof}

We finish with a generalization of a result of Iwanaga, see
\cite[Proposition 9.1.10]{Rel-hom}.

\begin{corollary}\label{geneIwana}
Assume that $\lGgldim(R)\leq n$ holds for some non-negative
integer $n$. If for an R-module M one of the numbers $\pd_R(M)$,
$\id_R(M)  $, or $\fd_R(M)$ is finite, then all of them are
smaller or equal to $n$.
\end{corollary}
\begin{proof}
If $\pd_R(M) $ is finite, then \cite[Proposition 2.27]{HH} and the
assumption give $\pd_R(M) = \Gpd_R(M)\leq n$. The argument for
$\id_R(M) < \infty$ is similar. Finally, Corollary
\ref{cor-main}(2) and the assumption give $\lFPD(R) \leq n$, and
then $\fd_R(M) < \infty$ implies  $\pd_R(M) < \infty$ by
\cite[Proposition 6]{Jensen}.\end{proof}

\noindent {\bf ACKNOWLEDGEMENTS.} The authors would like to express their sincere thanks for
the referee for his/her helpful suggestions. \\

\bibliographystyle{amsplain}

\end{document}